\documentclass[12pt]{article}

\usepackage{amsfonts}
\usepackage{amsmath,amssymb,amsfonts}
\usepackage{color}
\usepackage{hyperref}
\usepackage[latin1]{inputenc}
\usepackage{amsthm}
\usepackage{ulem}
\usepackage{xfrac}
\usepackage{tikz-cd}
\usepackage{graphicx}
\usepackage{caption}
\usepackage{subcaption}
\usepackage{float}

\newcommand{\mg}{\mathfrak g }

\newcommand{\ms}{\mathfrak s }

\newcommand{\mmu}{\mathfrak u }

\newcommand{\mk}{\mathfrak k }

\newcommand{\ma}{\mathfrak a }

\newcommand{\mpp}{\mathfrak p }

\newcommand{\so}{\mathfrak{so} }
\newcommand{\su}{\mathfrak{su} }

\renewcommand{\sl}{\mathfrak{sl} }

\newcommand{\gl}{\mathfrak{gl} }

\newcommand{\R}{\mathbb R}

\newcommand{\N}{\mathbb N}

\newcommand{\C}{\mathbb C}

\def\squarebox#1{\hbox to #1{\hfill\vbox to #1{\vfill}}}
\newcommand{\qedem}{\hspace*{\fill} \vbox{\hrule\hbox{\vrule\squarebox{.667em}\vrule}\hrule}\smallskip}
\newtheorem{teo}{Theorem}[section]
\newtheorem{prop}[teo]{Proposition}

\newtheorem{lema}[teo]{Lemma}
\newtheorem{defin}[teo]{Definition}
\newtheorem{explo}[teo]{Example}

\newenvironment{profe}{\noindent{\bf Proof:}}{\hfill\qedem\vspace{3mm}}

\begin{document}

\title{On invariant control sets for  control systems on  $S^3$}
\author{Bruno A. Rodrigues\\State University of Maringa, Brazil \and Luiz A. B. San Martin\\University of Campinas, Brazil \and  Alexandre J. Santana\\State University of Maringa, Brazil}
\maketitle
\textbf{Abstract:}  In this paper we describe the Lie-theoretic structure of ${\rm SO}(1,4)$ and consider control systems given by certain vector fields of  ${\rm SO}(1,4)$. Then we explicitly describe  its  invariant control sets in the unique ${\rm SO}(1,4)$-flag manifold, namely the sphere $S^3$.


\textbf{MSC\ 2020.} 93B05, 15B30, 57M60, 57R27.

\section{Introduction}

A control system, the main concept  of geometric  control theory, on a manifold $M$ is a family of vector fields on $M$. One of the primary goal of control theory is to  develop control strategies that allow to determine which regions of $M$ where the system is controllable.
Recall that a control system is controllable if it sends any state (point of $M$) of the system to any  other state  in finite time duration. Despite the large amount of papers  and books on controllability, this subject is far from saturation as most results are still partial (see e.g. Elliott \cite{Elliott}).

On the other hand, if we know that the system is not controllable we can try to study maximal regions of the state space $M$ in which interior the controllability property
holds. This region, called invariant control set, was introduced by Arnold and Kliemann in the context of stochastic control theory (see \cite{ArKli}) and then  applied in the context of control theory (see e.g. Colonius and Kliemann \cite{ColK00}) and in the  context of actions of semigroups of Lie groups (see e.g. San Martin \cite{flagtype} and references therein). But is not easy to find  literature  with  explicit descriptions of invariant control sets  of  control systems on flag manifolds, in general due to the difficulty of understanding the geometry and topology of flag manifolds and describing the behavior of vector fields and their integral curves. Although in recent years a line of research has been successful in describing control sets in a particular context, namely in case of  linear systems on Lie groups (see e.g.  Ayala and Da Silva \cite{AySilva} and references therein).

The main contribution of our paper is to explicitly describe  the  invariant control sets (in the sphere $S^3$) given from control systems of certain vector fields of  ${\rm SO}(1,4)$. To achieve this objective we also presented  a Lie-theoretic setup of the Lie group ${\rm SO}(1,4)$.
Our approach is the following. First, note that a Cartan decomposition of the Lie algebra $\so(1,4)$ can be identified as the direct sum $\so(4)\oplus\mathbb{H}$, in which $\so(4)$ corresponds to the compact component, and $\mathbb{H}$ stands for the symmetric one. The maximal abelian subalgebra contained in the symmetric component is one-dimensional, implying that there is only one flag manifold for ${\rm SO}(1,4)$, which is precisely the sphere $S^3$. For symmetric elements, the vector fields given by the infinitesimal action of $\so(1,4)$ on $S^3$ are gradient vector fields of height functions, and elements in the compact component $\so(4)$ give rise to vector fields defined by right and left multiplication by imaginary quaternions. From this setting, the quaternion algebra shows to be highly useful for exploring invariant control sets for vector fields induced on $S^3$ by ${\rm SO}(1,4)$. 

Precisely, we consider control systems with drifts being vector fields corresponding to the real quaternion ${\bf 1}\in\mathbb{H}$ and control vector fields corresponding to pure quaternions. The control sets for such systems appear as spherical domes in some cases, while in others, they are described as geodesically convex closures of the set of attractor points for the vector fields corresponding to the control system. This observation seems to be true in general for appropriate families of vector fields whose trajectories follow geodesics. In this direction, we establish that if the set $E$ of attractors for a family $\Gamma$ of geodesic vector fields on a differentiable manifold $M$ is both closed and geodesically convex, then $E$ is the invariant control set for $\Gamma$.

The paper is structured in the following way. In the Section 2 we give part of the notations e definitions necessary in the paper. In section 3 we prove a result that establish conditions on vector fields to obtain a description of the  invariant control set on a differentiable manifold. The fourth section is dedicated to present  the structure of the Lie group ${\rm SO}(1, 4)$, the ${\rm SO}(1, 4)$-flag manifold and the vector fields on this manifold.  Finally, in section 5 we express the shape of the control sets for a family of control systems on
the ${\rm SO}(1, 4)$-flag manifold.

\section{Background}

Recall that a control system is a family of a complete vector fields $\Gamma $ on an $n$-dimensional manifold $M$. A trajectory of $\Gamma $ is a continuous curve $\gamma $ from an interval $[0,T]$, $T \geq 0$ of the real line into $M$ such that for some partition $0<t_{1} < t_{2} < \cdots < t_{n}=T$ exist vector fields $X_{1}, \ldots , X_{n}$ in $\Gamma$ such that the restriction of $\gamma$ to each interval $[t_{i-1}, t_{i} )$ is an integral curve of $X_{i}$. We consider in this paper a special case given by the equations 
\[ \dot{x}(t)=X_{1}(x)+ \Sigma_{i=2}^{l} u_{i}X_{i}(x), \]
where $X_{j} \in \Gamma$ and
$	\mathcal{U}=\{u\in L^{\infty}(\mathbb{R},\mathbb{R}^{m})\left\vert\, u(t)\in
	U\text{ for almost all }t\right.  \}$. As usual,  $X_{1}$ is called drift and $X_{2}, \ldots , X_{l}$ control vector fields.

 We assume that the control range $U \subset {\mathbb R}^{m}$ is nonvoid and that for every initial state $x \in M$ and every $u \in \mathcal{U}$ there exists a unique solution denoted by $\phi (t,x,u), t \in {\mathbb R}$, satisfying $\phi (0,x,u)=x$. We can also admit piece-wise constant controls, that is,
 \[
	\mathcal{U}_{pc}=\{u:\mathbb{R}\to U\subset\mathbb{R}^{m}\left\vert\, u\mbox{  piece-wise constant}\right.  \}.
	\]

In particular, we  are concerned here with  systems where the state space $M$ is a Lie group $G$ and $X_{i}$ are  (right or left) invariant vector fields on $G$ and with systems that the state space is a flag manifolds of $G$.
\begin{defin}
An invariant control system $\Gamma\subset\mg$ is said to satisfy the Lie algebra rank condition (LARC) if the Lie algebra generated by $\Gamma$, ${\rm Lie}(\Gamma)$, is the whole $\mg$, that is,
$${\rm Lie}(\Gamma)=\mg.$$
\end{defin}

The positive orbit of $x\in M$ at time exactly $t>0$ is the set
\[
\mathcal{O}^{+}_{t}(x) =\left\{ y \in M\,|\,\exists u\in\mathcal{U} \mbox{ with } y=\phi(t,x,u) 		  \right\}.
\]
Similarly, the negative orbit of $x\in M$ at time $t>0$ is
\[
\mathcal{O}^{-}_{t}(x) =\left\{ y \in M\,|\,\exists u\in\mathcal{U} \mbox{ with } x=\phi(t,y,u) 		  \right\}.
\]
The positive and negative orbits of $x\in M$ up to time $T$ are defined as
\begin{eqnarray*}
\mathcal{O}^{+}_{\leq T}(x) &=&\bigcup\nolimits_{0\leq t\leq T}\mathcal{O}^{+}_{t}(x)\\
&=&\left\{  y \in M \mbox{ such that }\exists\, 0 \leq t \leq T \mbox{ and }\exists u\in\mathcal{U} \mbox{ with } y=\phi(t,x,u) 		  \right\}
\end{eqnarray*}
and
\begin{eqnarray*}
\mathcal{O}^{-}_{\leq T}(x)&=&\bigcup\nolimits_{0\leq t\leq T}\mathcal{O}^{-}_{t}(x)\\
&=&\left\{ y \in M \mbox{ such that }\exists\, 0 \leq t \leq T \mbox{ and }\exists u\in\mathcal{U} \mbox{ with } x=\phi(t,y,u) 		  \right\}
\end{eqnarray*}
The positive and negative orbits of $x\in M$ are 
\[
\mathcal{O}^{+}(x) =\bigcup\nolimits_{T>0}\mathcal{O}^{+}_{\leq T}(x)=\bigcup\nolimits_{t>0}\mathcal{O}^{+}_{t}(x)
\]
and
\[
\mathcal{O}^{-}(x) =\bigcup\nolimits_{T>0}\mathcal{O}^{-}_{\leq T}(x)=\bigcup\nolimits_{t>0}\mathcal{O}^{-}_{t}(x).
\]

One of the key concepts of this work is controllability of  control system $\Gamma$, which roughly speaking signifies that the orbit of every single point of $M$ under $\Gamma$ covers the whole manifold $M$.
\begin{defin}
A control system $\Gamma$ is said to be controllable from $x\in M$ when $\mathcal{O}^+(x)=M$, and it is said to be controllable when it is controllable from every $x\in M$.
\end{defin}
When a control system $\Gamma$ fails to be controllable, one can ask for the maximal subsets of MM where controllability holds. These are the control sets. Complete controllability can occur only on connected Lie groups and the LARC   is a necessary condition for controllability, although in general it is not sufficient. 
\begin{defin}
\label{Definition_invcontrol_sets}

A nonvoid set $C \subset$ $M$ is an invariant control set of a control system $\Gamma$ on $M$ if it has the following properties: 
  \begin{itemize}
    \item[(i)] for all $x \in C$ there is a control $u\in\mathcal{U}$ such that $\phi(t,x,u)\in C$ for all $t\geq0$,
    \item[(ii)] for all $x\in C$ one has ${\rm cl}C = {\rm cl}(\mathcal{O}^{+}(x))$, and 
    \item[(iii)] $C$ is maximal with these properties, that is, if $C^{\prime}\supset C$ satisfies conditions (i) and (ii), then $C^{\prime}=C$.
  \end{itemize}
  \end{defin}

\section{General Theorem}

A subset $C$ of a differentiable manifold $M$ is geodesically convex if the minimal geodesic segment joining $p,q\in C$ is contained in $C$. We define also the geodesic convex hull of a subset as the smaller geodesically convex subset $C\subset M$ containing it.

Let $C^{\infty}(M)$ denote the set of all complete $\mathcal{C}^{\infty}$ vector fields on a differentiable manifold $M$. Given $X\in C^{\infty}(M)$, we say that a point $a\in M$ is an attractor fixed point  for $X$ in $M$ if the vector field $X$ vanishes in $a$ and $\lim_{t\to\infty}X_t(x)=a$ for every initial state $x\in M$, where $X_t$ denotes the flow corresponding to $X$. For simplicity we will refer to such points as just the attractors of $X$.

\begin{teo}
Let $M$ be a differentiable manifold with an affine connection $\nabla$ and $\Gamma\subset C^{\infty}(M)$ a family of vector fields on $M$ such that every $X\in\Gamma$ satisfies $\nabla_XX=c_XX$, where $c_X:M\to\R$ is a smooth function on $M$. Assume that every $X\in\Gamma$ admits only one  attractor  $a\in M$ and suppose also that a trajectory $\gamma$ for $X$ is minimal geodesic between all of its points and the attractor $a$. If
$$E=\left\{a\in M\,|\, a\hbox{ is attractor  for some }X\in\Gamma\right\}$$
is closed and geodesically convex, then $E$ is the invariant control set for $\Gamma$.
\label{geral}
\end{teo}

\begin{profe}
The condition $\nabla_XX=c_XX$ on the vector fields of the system tells us that the trajectories for the vector fields $X\in\Gamma$ follow the geodesics on $M$, and in this case we call $X$ a geodesic vector field. 

Now, let $D$ be the invariant control set for $\Gamma$. We know that 
$$D=\bigcap_{x\in M}{\rm cl}(S_{\Gamma}\cdot x),$$
where $S_{\Gamma}$ denotes the semigoup of the system $\Gamma$. Since every $y\in E$ is an attractor  for some $X\in\Gamma$, we have $y\in{\rm cl}(S_{\Gamma}\cdot x)$, for every $x\in M$. In other words, $E\subseteq D$. Being $D$ a control set, if $x\in D$ then $x\in{\rm cl}(S_{\Gamma}\cdot y)$, for every $y\in D$. In particular, $x\in{\rm cl}(S_{\Gamma}\cdot y)$ for every $y\in E$ (we just proved that $E\subseteq D$). Well, $x\in{\rm cl}(S_{\Gamma}\cdot y)$ implies that we can choose a trajectory $\gamma$ for the control system $\Gamma$ with starting point at $y$ and endpoint $x_T$ arbitrarily close from $x$. Let us say that the very last path $\gamma_T$ of $\gamma$ joins $x_{T'}\in E$ to $x_T$ and corresponds to $X_0\in\Gamma$ with attractor  $x_0$. By the definition of $E$ we have that $x_0\in E$ and $\gamma_T$ is a minimal geodesic segment joining $x_{T'}$ to $x_0$. Since $A$ is geodesically convex, we must have $\gamma_T\subset E$, that is $x_T\in E$. That is, $x\in{\rm cl}E$, and this proves the equality $E=D$. Hence $E$ is the invariant control set for $\Gamma$.
\end{profe}

The above theorem assumes the set of attractors $E$ to be geodesically convex. Nonetheless, there are examples where this convexity does not hold.

In the course of this work we present some examples, which suggest that if this occurs then the convex hull of $E$ can be a good candidate for the invariant control set. Our first example will consider a control system evolving on the projective space $\mathbb{P}^{n-1}$ induced by the action of ${\rm Sl}(n,\R)$.

\begin{explo}

Let 
$$v_0=\left(\begin{array}{c}\frac{1}{\sqrt{n}}\\ w\end{array}\right)\in\R^n\hbox{ with }\frac{1}{n}+|w|^2=1$$ 
and 
$$B_i=\left[\begin{array}{cc}0&0\\0&X_i\end{array}\right]\in\sl(n,\R),\hspace{5mm}i=1,\dots,m=\dim(\so(n-1)),$$
where each diagonal block $X_i$ has order $n-1$ and $\{X_1,\dots,X_m\}$ is a basis for $\so(n-1)$. Setting $A=v_0v_0^t-\frac{1}{n}{\rm Id_n}$, we will describe the set of attractors $E$ and the invariant control set $C$ for the control system 
\begin{equation}
\dot{x}=Ax+u_1B_1x+\dots+u_mB_mx,\hspace{5mm}x\in\mathbb{P}^{n-1}.
\label{sistema}
\end{equation}

Before anything else, note that (\ref{sistema}) satisfies the LARC. In fact, to verify this let ${\rm ad}$ be the adjoint representation  of
$$\so(n-1)=\left[\begin{array}{cc}0&0\\0&\so(n-1)\end{array}\right]$$
on $\ms\subset\sl(n,\R)$, the subspace of symmetric matrices.
The restriction of ${\rm ad}$ to the subspace
$$W=\left\{\left[\begin{array}{cc}0&w^t\\ w&0\end{array}\right];w\in\mathbb{R}^{n-1}\right\}$$
of $\sl(n,\R)$ is an irreducible representation of ${\rm ad}$. In fact, if
$$B=\left[\begin{array}{cc}0&0\\0&X\end{array}\right]\in\so(n-1)\hbox{ and }\tilde{w}=\left[\begin{array}{cc}0&w^t\\w&0\end{array}\right]$$
we have
$$[B,\tilde{w}]=\left[\begin{array}{cc}0&(Xw)^t\\ Xw&0\end{array}\right].$$
Similarly, the restriction of ${\rm ad}$ to 
$$\ms_0=\left\{\left[\begin{array}{cc}0&0\\0&S\end{array}\right];\,\,S\hbox{ traceless symmetric}\right\}$$
is irreducible. Now, we can decompose $A=v_0v_0^t-\frac{1}{n}{\rm Id}_n$ as
$$A=\frac{1}{\sqrt{n}}\left[\begin{array}{cc}0&w^t\\w&0\end{array}\right]+\left[\begin{array}{cc}0&0\\0&ww^t-\frac{1}{n}{\rm Id}_{n-1}\end{array}\right]=:\frac{1}{\sqrt{n}}\tilde{w}+S_0,$$
and since $w\neq0$ we obtain that the subspaces generated by
$${\rm ad}\left(\so(n-1)\right)\left(\frac{1}{\sqrt{n}}\tilde{w}\right)\hbox{ and }{\rm ad}\left(\so(n-1)\right)\left(S_0\right)$$
are $W$ and $\ms_0$, respectively. Finally, the skew-symmetric matrices of the form
$$\left[\begin{array}{cc}0&-w^t\\w&0\end{array}\right],\hspace{3mm}w\in\R^{n-1},$$
are obtained as Lie brackets between suitable matrices 
$$\left[\begin{array}{cc}0&u^t\\u&0\end{array}\right]\in\ms_0\hbox{ and }\left[\begin{array}{cc}0&0\\0&X\end{array}\right]\in\so(n-1),$$
proving that $(\ref{sistema})$ satisfies LARC, and this implies that the invariant control set for this system has nonempty interior.

Now we describe the system (\ref{sistema}) in a more detailed way. For, let ${\rm Gl}(n,\R)$ acting on the projective space $\mathbb{P}^{n-1}$ under $g*[x]:=[gx]$, where $g\in{\rm Gl}(n,\R)$ and $[x]\in\mathbb{P}^{n-1}$ indicates the straight line containing $0$ and $x\in S^{n-1}$. We have that 
$$g*x=\left[\frac{gx}{|gx|}\right],$$ 
and given any $X\in\gl(n,\R)$ the corresponding infinitesimal action on $\mathbb{P}^{n-1}$ is
$$\tilde{X}[x]=\frac{d}{dt}\left.\left(e^{tX}*x\right)\right|_{t=0}=\frac{d}{dt}\left.\left(\frac{e^{tX}x}{|e^{tX}x|}\right)\right|_{t=0}=[Xx-\langle Xx,x\rangle x].$$

For instance, if $X={\rm diag}(1,0,\dots,0)\in\gl(n,\R)$ and $x=(x_1,\dots,x_n)\in S^{n-1}$, then $Xx=(x_1,0,\dots,0)$ is the projection of $x$ along $e_1=(1,0,\dots,0)$ and $\tilde{X}[x]$ is given by the orthogonal projection of $Xx$ on the one dimensional subspace generated by $x$. We have that $\tilde{X}[e_1]=0$ and if $x$ belongs to the orthogonal complement of $e_1$ we also have $\tilde{X}[x]=0$. The trajectories for $\tilde{X}$ follow the great circles and $[e_1]$ is the only attractor fixed point for $\tilde{X}$ while the points in the orthogonal complement of $e_1$ are the reppeler fixed points for $\tilde{X}$.

In this way, if $\Gamma\subset\mathbb{P}^{n-1}$ is a closed and geodesically convex set, for every $[x]\in\Gamma$ corresponds a vector field $\tilde{X}_x$ where $X_x$ indicates the projection over the line generated by $x$ in $\R^n$. If $\tilde{\Gamma}$ is the family of such vector fields, then the previous theorem tells us that its invariant control system is exactly $\Gamma$.

Now, let $\ms\subset\sl(n,\R)$ the subspace of symmetric matrices, and consider the embedding of the projective space $\mathbb{P}^{n-1}$ in $\ms$ which associates to each $[v]\in\mathbb{P}^{n-1}$, $|v|=1$, the symmetric $n\times n$ trace free matrix $V$ given by
$$V=vv^t-\frac{1}{n}{\rm Id}_n.$$

Given $[v],[w]\in\mathbb{P}^{n-1}$, with respect to the canonical inner product in $\ms$, we have
\begin{eqnarray*}
\langle V,W\rangle&=&{\rm tr} VW\\
&=&{\rm tr}\left(vv^tww^t-\frac{1}{n}vv^t-\frac{1}{n}ww^t+\frac{1}{n^2}{\rm Id}_n\right)\\
&=&{\rm tr}\left(\langle v,w\rangle vw^t-\frac{1}{n}vv^t-\frac{1}{n}ww^t+\frac{1}{n^2}{\rm Id}_n\right)\\
\end{eqnarray*}
since ${\rm tr}(vv^t)=|v|$ and ${\rm tr}(vw^t)=\langle v,w\rangle$, we get
$$\langle V,W\rangle=\langle v,w\rangle^2-\frac{1}{n}-\frac{1}{n}+\frac{n}{n^2}=\langle v,w\rangle^2-\frac{1}{n}.$$

For the symmetric matrix $E^0_{11}=e_1e_1^t-\frac{1}{n}{\rm Id}_n$, let the $C$ the convex set in $\mathbb{P}^{n-1}$ defined by
$$C=\left\{[v]\in\mathbb{P}^{n-1}\,|\,\langle V,E^0_{11}\rangle\geq0\right\}.$$
We have
$$C=\left\{[v]\in\mathbb{P}^{n-1}\,|\,\langle v,e_1\rangle^2-\frac{1}{n}\geq0\right\}=\left\{[v]\in\mathbb{P}^{n-1}\,|\,\langle v,e_1\rangle\geq\frac{1}{\sqrt{n}}\right\},$$
and if $\theta$ is the angle between $v$ and $e_1$, then $C$ is the set of $[v]\in\mathbb{P}^{n-1}$ such that $\cos(\theta)\geq\sfrac{1}{\sqrt{n}}$, that is, the spherical dome in $\mathbb{P}^{n-1}$ of the elements $[v]$ whose $v_1\geq\sfrac{1}{\sqrt{n}}$, and $C$ can be seen as the intersection of an halfspace in $\ms$ with $\mathbb{P}^{n-1}$ (using the embedding previously defined).

Returning to system (\ref{sistema}), note that $[v_0]\in\mathbb{P}^{n-1}$ is the attractor fixed point for the vector field on $\mathbb{P}^{n-1}$ induced by $A=v_0v_0^t-\frac{1}{n}{\rm Id}$. In fact, $|v_0|=1$ and
\begin{eqnarray*}
Av_0-\langle Av_0,v_0\rangle v_0&=&\left(v_0v_0^t-\frac{1}{n}{\rm Id}_n\right)v_0-\left\langle\left(v_0v_0^t-\frac{1}{n}{\rm Id}_n\right)v_0,v_0\right\rangle v_0\\
&=&v_0v_0^tv_0-\frac{1}{n}v_0-\left\langle v_0v_0^tv_0-\frac{1}{n}v_0,v_0\right\rangle v_0\\
&=&|v_0|v_0-\frac{1}{n}v_0-\left(|v_0|^3-\frac{1}{n}|v_0|^2\right)v_0\\
&=&0
\end{eqnarray*}
Also, the spherical dome $C$ is invariant under each $B_i$, since
$$e^{tB_i}=\left[\begin{array}{cc}1&0\\0&e^{tX_i}\end{array}\right]$$
and
$$\left\langle e^{tB_i}v,e^{tB_i}e_1\right\rangle=\left\langle e^{tB_i}v,e_1\right\rangle\geq\frac{1}{\sqrt{n}},$$
for every $v\in\R^n$ such that $\langle v,e_1\rangle\geq\sfrac{1}{\sqrt{n}}$.

Being $\{X_1,\dots,X_m\}$ a basis for $\so(n-1)$, choosing adequate controls we see that the attractors for the vector fields induced by $A+u_1B_1+\cdots+u_mB_m$ in $\mathbb{P}^{n-1}$ just rotate around the boundary of $C$, meaning that
$$E=\left\{[v]\in\mathbb{P}^{n-1}\,|\,\langle v,e_1\rangle=\frac{1}{\sqrt{n}}\right\}$$
is fulfilled with attractors for system (\ref{sistema}), that is, $C$ is invariant under $(\ref{sistema})$, and hence the invariant control set for this system must be contained in $C$. Furthermore, if $\mathcal{O}^+(x)$ stands for the positive orbit of $x\in\mathbb{P}^{n-1}$, then for every $y\in E$ we have 
$$y\in\bigcap_{x\in \mathbb{P}^{n-1}}{\rm cl}\mathcal{O}^+(X),$$
ensuring that $C$ is the invariant control set for $(\ref{sistema})$.

\end{explo}

\section{Lie-theoretic structure of \texorpdfstring{${\rm SO}(1,4)$}{SO(1,4)}}

The purpose of this section is to present the concepts and structures necessary to describe the control systems for which we will study the invariant control sets. As it is difficult to find bibliographical references that deal with Lie theory in the special case of ${\rm SO}(1,4)$, then we will provide here the necessary theory of ${\rm SO}(1,4)$ and of its Lie algebra for the purpose of this work, for example, the description of the Cartan decomposition, the flag manifold and the vector fields on this manifold.

\subsection{Cartan decomposition}

The quadratic form with matrix 
$$I_{p,q}=\left[\begin{array}{cc}1_p&0\\0&-1_q\end{array}\right]$$
gives rise to the indefinite special orthogonal Lie algebra of type $(p,q)$. In other words, the Lie algebra $\so(p,q)$ consists of the order $(p+q)$ real matrices that are skew-symmetric with respect to $I_{p,q}$, that is, 
$$\so(p,q)=\{X\in\sl(p+q,\R)\,|\,I_{p,q}X+X^tI_{p,q}=0\}.$$ 
The skew-symmetry with respect to any other quadratic form equivalent to $I_{p,q}$ defines a Lie algebra isomorphic to $\so(p,q)$.
We are interested in the Lie algebra $\so(1,4)$, and in this case the condition $I_{1,4}X+X^tI_{1,4}=0$ implies that the matrix $X$ (which is an order $5$ matrix) must have the block form  
$$X=\left[\begin{array}{cc}0&\beta\\ \beta^t&\gamma\end{array}\right],\hbox{ with }\gamma=-\gamma^t,\hbox{ that is, }\gamma\in\so(4)=\su(2)\oplus\su(2).$$
The corresponding compact real form in this case is the Lie algebra $\so(5)$ and the associate conjugation in $\so(5,\C)$ is $\theta(X)=-X^t$ (note that it is an involutive automorphism for this Lie algebra). The bilinear form $B_{\theta}$ in $\so(1,4)$ given by
$$B_{\theta}(X,Y)=-\langle X,\theta Y\rangle$$
is an inner product in $\so(1,4)$, as follows from $\cite{SMALG}$ (Lemma 12.21). This implies that $\so(1,4)$ decomposes as the direct sum of the eigenspaces
$$V_1=\{X\in\so(1,4)\,|\,\theta(X)=X\}\hbox{ and }V_{-1}=\{X\in\so(1,4)\,|\,\theta(X)=-X\}.$$
Setting $\mk=V_1$ and $\ms=V_{-1}$, then $\so(1,4)=\mk\oplus\ms$ is a Cartan decomposition. Explicitly, we have
$$\mk=\left\{\left.\left[\begin{array}{cc}0&0\\0&\gamma\end{array}\right]\in\so(1,4)\right|\gamma\in\so(4)\right\}=\so(4)$$
and
$$\ms=\left\{\left.\left[\begin{array}{cc}0&\beta\\ \beta^t&0\end{array}\right]\in\so(1,4)\right|\beta\in\R^4\right\}.$$
Also, $\mmu=\mk\oplus{\bf i}\ms$ is a compact real form for $\so(5,\C)$ (see $\cite{SMALG}$ Proposition 12.27). It is given by
$$\mmu=\left\{\left.\left[\begin{array}{cc}0&{\bf i}\beta\\ {\bf i}\beta^t&\gamma\end{array}\right]\right|\beta\in\R^4,\,\,\gamma\in\so(4)\right\}.$$

The application $\mmu\to\so(5)$ defined by 
$$\left[\begin{array}{cc}0&{\bf i}\beta\\ {\bf i}\beta^t&\gamma\end{array}\right]\mapsto\left[\begin{array}{cc}0&\beta\\ -\beta^t&\gamma\end{array}\right]$$
is a Lie isomorphism between $\mmu$ and $\so(5)$.
A maximal abelian subalgebra $\ma$ contained in $\ms$ is generated by the matrices
$$\left[\begin{array}{cc}0&\beta\\ \beta^t&0\end{array}\right],\hspace{2mm}\hbox{where }\beta=[1,0,0,0],$$
and thus $\so(1,4)$ is a real rank 1 non-compact real form of $\so(5,\C)$.

When realizing $\so(1,4)$ by the quadratic form given by 
$$J_{1,n}=\left(\begin{array}{ccc}-1_{3}&0&0\\ 0&0&1\\ 0&1&0\end{array}\right)$$
we get a simpler way to compute the restricted roots with respect to a maximal abelian subalgebra in $\ms$.
In this realization $X\in\so(1,4)$ if, and only if, $X$ has the form 
$$\left(\begin{array}{ccc}A&B&C\\ C^t&\alpha&0\\ B^t&0&-\alpha\end{array}\right),\,\,\hbox{ where }A\in\so(3).$$
A Cartan decomposition is given by the symmetric and the anti-symmetric matrices of this type. Precisely,
$$\mk=\left\{\left.\left[\begin{array}{ccc}A&B&-B \\ -B^t&0&0\\ B^t&0&0\end{array}\right]\,\right|\, A\in\so(3),\,B^t\in\R^3\right\}$$
and
$$\ms=\left\{\left.\left[\begin{array}{ccc}0&B&B \\ B^t&\alpha&0\\ B^t&0&-\alpha\end{array}\right]\,\right|\, B^t\in\R^3,\,\alpha\in\R\right\}.$$

In this realization, a maximal abelian subalgebra contained in $\ms$ is
$$\ma=\left\{\left.\left[\begin{array}{ccc}0&0&0 \\ 0&\alpha&0\\ 0&0&-\alpha\end{array}\right]\,\right|\, \alpha\in\R\right\}.$$
And the restricted roots are the functionals $\pm\lambda:\ma\to\R$, $\lambda(H)=\alpha$, where $H={\rm diag}(0,\alpha,-\alpha)\in\ma.$
The only simple root is $\lambda$ and it has multiplicity 3. The corresponding root spaces are
\begin{eqnarray*}
\mg_{\lambda}&=&\{X\in\so(1,4)\,|\,{\rm ad}(H)X=\lambda(X),\,\,\forall H\in\ma\}\\
&=&\left\{\left.\left[\begin{array}{ccc}0&0&C\\C^t&0&0\\0&0&0\end{array}\right]\in\so(1,4)\,\,\right|\,\,C^t\in\R^3\right\},
\end{eqnarray*}
similarly we have
\begin{eqnarray*}
\mg_{-\lambda}&=&\left\{\left.\left[\begin{array}{ccc}0&B&0\\0&0&0\\B^t&0&0\end{array}\right]\in\so(1,4)\,\,\right|\,\,B^t\in\R^3\right\}
\end{eqnarray*}
and
\begin{eqnarray*}
\mg_{0}&=&\left\{\left.\left[\begin{array}{ccc}A&0&0\\0&\alpha&0\\0&0&-\alpha\end{array}\right]\in\so(1,4)\,\,\right|\,\,A\in\so(3),\,\,\alpha\in\R\right\}.
\end{eqnarray*}

\subsection{The flag manifold}

Flag manifolds are connected (because $G$ is connected) and the following proposition ensures that flag manifolds are compact as well.

\begin{prop} \label{Kact}
$K$ acts transitively on any flag manifold $\mathbb{F}_{\Theta}$.
\end{prop} 
\begin{profe}
First of all, we know that $\mathbb{F}_{\Theta}=G/P_{\Theta}$ has Lie algebra $\mg/\mpp_{\Theta}$, and $\dim\mathbb{F}_{\Theta}=\dim{\mg}-\dim{\mpp_{\Theta}}$. Let $\pi:G\to G/P_{\Theta}$ be the canonical projection. The orbit $K\cdot b_0$ of the origin $b_0\in G/P_{\Theta}$ is a submanifold of $G/P_{\Theta}$ with tangent space $d\pi_1(\mk)$, and $\dim d\pi_1(\mk)=\dim\mk-\dim(\mk\cap\mpp)$. Now, the Iwasawa decomposition gives us $\mg=\mk+\mpp_{\Theta}$, hence $\dim\mg=\dim\mk+\dim\mpp_{\Theta}-\dim(\mk\cap\mpp_{\Theta})$. This implies that the dimension of the orbit $K\cdot b_0$ coincides with $\dim(G/P_{\Theta})$, and hence the orbit is an open submanifold in $G/P_{\Theta}$. Since $K$ is compact, we get that $K\cdot b_0$ is closed as well. Finally, the connectedness of $\mathbb{F}_{\Theta}=G/P_{\Theta}$ implies that $K\cdot b_0=\mathbb{F}_{\Theta}$, that is, $K$ acts transitively on the flag manifold $\mathbb{F}_{\Theta}$.  
\end{profe}

By Proposition \ref{Kact} we have that $\mathbb{F}_{\Theta}=K/K_{\Theta}$, where the stabilizer is given by $K_{\Theta}=K\cap P_{\Theta}$.

Since $\so(1,4)$ is a rank one real Lie algebra, there is only one flag manifold for ${\rm SO}(1,4)$, namely the maximal one. Note that we can identify $\mathbb{F}=\sfrac{SO(1,4)}{P}$ as a $K$-orbit under the adjoint representation, that is, ${\rm Ad}(K)H=\sfrac{K}{K_H}$, where $H\in{\rm cl}\ma^+$ and $K_H$ is the centralizer of $H$ in $K$. 

\begin{prop}
The only flag manifold $\mathbb{F}$ of $SO(1,4)$ embeds in the component $\ms$ of the Cartan decomposition as the ${\rm Ad}(SO(4))$-orbit of $H\in\ma^+$.
\end{prop}
\begin{profe}
The component $\ms$ of the Cartan decomposition $\mg=\mk\oplus\ms=\so(4)\oplus\ms$ is invariant under the adjoint representation of $K= SO(4)$. Since $\ma$ is a one dimensional subalgebra, a Weyl chamber $\ma^+\subset\ma$ is just a ray starting at the origin. Choosing an element $H\in\ma^+$, the stabilizer of $H$ under the adjoint action of $K$ on $\ms$ is the centralizer $K_H$ of $H$ in $K$, which is given by $K\cap P$. It follows that the adjoint orbit ${\rm Ad}(SO(4))H\subset\ms$ is identified with the coset space $K/K_H$, and this one is the  flag manifold $\mathbb{F}=G/P$.
\end{profe}

\begin{prop}
The sphere $S^3$ is the only flag manifold of ${\rm SO}(1,4)$.
\end{prop}
\begin{profe}
We look at $\so(1,4)$ realized by $I_{1,4}$, where the maximal abelian subalgebra $\ma$ contained in $\ms$ is of the form 
$$\ma= {\rm span}\left\{ \left[\begin{array}{cc}0&e_1\\e^t_1&0\end{array}\right]  ,\,\,\,e_1=(1,0,0) \right\} .$$
           
Let 
$$H=\left[\begin{array}{cc}0&e_1\\e_1^t&0\end{array}\right]\in\ma^+\,\,\,\hbox{ and }\,\,\,x=\left[\begin{array}{cc}1&0\\0&\beta\end{array}\right]\in K, \mbox{ with } \beta\in{\rm SO}(4).$$
Denote by  $\beta_{1i}$ and by $\beta_{i1}$  the first row and the first column of $\beta$ respectively.  With this notation we have
$$xH=Hx\Leftrightarrow \left[\begin{array}{cc}0&e_1\\\beta_{i1}&0\end{array}\right]=\left[\begin{array}{cc}0&\beta_{1i}\\e_1^t&0\end{array}\right],$$
from which $\beta_{1i}=e_1$ and $\beta_{i1}=e_1^t$. This means that $\beta$ is of the form
$$\beta=\left[\begin{array}{cc}1&0\\0&\gamma\end{array}\right],\,\,\,\gamma\in{\rm SO}(3).$$
Hence, $x\in K_H$ if, and only if, $x$ has the block form
$$x=\left[\begin{array}{ccc}1&0&0\\0&1&0\\0&0&\gamma\end{array}\right],\,\,\,\gamma\in{\rm SO}(3).$$
This shows that $K_H={\rm SO}(3)$, and we get $\mathbb{F}=\sfrac{K}{K_H}=\sfrac{{\rm SO}(4)}{{\rm SO}(3)}= S^3$.
\end{profe}

The previous computations remain basically unchanged for ${\rm SO}(1,n)$, and in this case we have $\mathbb{F}=\sfrac{{\rm SO}(n)}{{\rm SO}(n-1)}= S^{n-1}$.

\subsection{Vector fields  on \texorpdfstring{$S^3$}{S3}}

Now we describe the vector fields induced on the sphere $S^3$ by the action of ${\rm SO}(1,4)$. Remember that an infinitesimal action of $\so(1,4)$ on $S^3$ is a homomorphism $\so(1,4)\to\Gamma(TS^3)$, where $\Gamma(TS^3)$ stands for the Lie algebra of vector fields on $S^3$. By means of an infinitesimal action of $\so(1,4)$ on $S^3$ one can see the Lie algebra $\so(1,4)$ as a Lie algebra of vector fields on the sphere $S^3$. 

\begin{teo}
The infinitesimal action of $\so(1,4)$ induced by the action of ${\rm SO}(1,4)$ on $S^3$ has as image the vector space formed by the vector fields
$$X_{(q,z,w)}(x)=\frac{1}{2}(q-x\overline{q}x)+zx+xw,\hspace{5mm}x\in S^3,$$
where $q\in\mathbb{H}=\ms$ and $z,w\in{\rm Im}\mathbb{H}=\su(2)$. This vector space is a $10$-dimensional Lie algebra  isomorphic to $\so(1,4)$.
\label{infiniact}
\end{teo}
\begin{profe}
We begin by investigating the vector fields corresponding to elements belonging to the $\mathfrak{s}$ component. There exists a $K$-invariant Riemannian metric such that for every $q \in {\mathbb H}=\mathfrak{s}$ the vector field $\tilde{X}_q$ induced by $q$ on $S^3$ is the gradient of the height function $f_q(\cdot)=\langle q,\cdot\rangle$ with respect to this $K$-invariant metric (see Duistermaat et al \cite{DKV} and Takeuchi and Kobayashi \cite{TK} for details). In the present case, since $\alpha(H)=1$ ($H\in{\mathfrak{s}}^+$) for every positive root $\alpha$ for which $\alpha(H)\neq0$ we have that the Borel metric coincides with the metric induced by the immersion of $S^3$ in $\mathfrak{s}$ (by this reason $S^3$ is called an immersed flag manifold).
The height function $f_q$ is linear on $\ms$, so its gradient vector field evaluated at $p\in S^3$ is obtained from the orthogonal projection of $q$ over $p$. In fact, $({\rm grad}f_q)_p=d(f_q)_p(v)$ is the cotangent vector $\omega$ such that $\omega(v)=\langle q,v\rangle$, that is, the cotangent vector $\omega$ such that 
$$\langle\omega,v\rangle=\langle q,v\rangle,\,\,\,\forall v\in T_pS^3.$$
Since
$$\langle q-\langle q,p\rangle p,v\rangle=\langle q,v\rangle\,\,\,\forall v\in T_pS^3,$$
we get 
$$({\rm grad}f_q)_p=q-\langle q,p\rangle p.$$
The vector field $\tilde{X}_q$ is thus given by
\begin{eqnarray*}
\tilde{X}_q(p)&=&q-\langle q,p\rangle p\\
&=&q-\frac{1}{2}\left(q\overline{p}+p\overline{q}\right)p\\
&=&q-\frac{1}{2}(q\overline{p}p+p\overline{q}p),\hspace{1cm}p\in S^3.
\end{eqnarray*}
Since $|p|=p\overline{p}=1$, we get
$$\tilde{X}_q(p)=\frac{1}{2}(q-p\overline{q}p).$$
Now we turn our attention to the elements in the compact component $\mk$. We know that the adjoint representation 
\begin{eqnarray*}
{\rm Ad}:K&\to&{\rm Gl}(\so(1,4))\\
k&\mapsto &{\rm Ad}(k)=d(C_k)_1.
\end{eqnarray*}
is differentiable and defines the action 
$$K\times\ms\to\ms,\,\,\, (k,x)\mapsto{\rm Ad}(k)x.$$
Since the flag $S^3$ embeds in $\ms$ as an ${\rm Ad}(K)$-orbit, we can consider the restriction of the above action to $S^3$ (viewed as an ${\rm Ad}(K)$-orbit) to get an infinitesimal action of $\mk$ on $S^3$. Thus, for $X\in\mk$ the corresponding infinitesimal action on $S^3$ is given by
$$\tilde{X}(x)=d({\rm Ad})_1(X)x={\rm ad}(X)x,$$
that is, the induced vector field is given by the adjoint action of $\mk$ on $\ms$. 

Now, as the Lie algebra ${\rm Im}\mathbb{H}$ is represented in $\mathbb{H}$ through left multiplication and also through right multiplication we have that the Lie algebra $\so(4)=\su(2)\oplus\su(2)$ is isomorphic to the Lie algebra of linear transformations
$$\{E_z+D_{w}\,|\,z,w\in{\rm Im}\mathbb{H}\}.$$

As the Lie algebra $\so(4)$ decomposes as the sum of two simple ideals commuting to each other, we describe the adjoint representation of $\mk=\so(4)$ on $\ms$ looking at each simple ideal. The first component corresponds to left multiplication by immaginary quaternions and the second one corresponds to right multiplication by immaginary quaternions, since these two kinds of quaternionic multiplications commute with each other as we have already seen.

We can check this as follows. Given $\beta=(p,q,r,s)\in\R^4$, we identify an element $S\in\ms$ as a quaternion number in the following way
\begin{equation}
S=\left[\begin{array}{cc}0&\beta\\\beta^t&0\end{array}\right]\in\ms\Leftrightarrow S=p+q{\bf i}+r{\bf j}+s{\bf k}\in\mathbb{H}.
\label{ident}
\end{equation}
Set
$$A_2=\left[\begin{array}{cc}0&-1\\1&0\end{array}\right],\,\,\,B_2=\left[\begin{array}{cc}0&1\\1&0\end{array}\right]\,\hbox{ and }\,C_2=\left[\begin{array}{cc}-1&0\\0&1\end{array}\right].$$
With this notation, define the $\so(4)$ matrices 
$$\begin{tabular}{ll}\vspace{5mm}
$\gamma_{\bf i}=\left[\begin{array}{cc}A_2&0\\0&-A_2\end{array}\right]$,   & $_{\bf i}\gamma=\left[\begin{array}{cc}A_2&0\\0&A_2\end{array}\right]$, \\ \vspace{5mm}
$\gamma_{\bf j}=\left[\begin{array}{cc}0&-1_2\\1_2&0\end{array}\right]$,   &  $_{\bf j}\gamma=\left[\begin{array}{cc}0&C_2\\-C_2&0\end{array}\right]$,  \\
$\gamma_{\bf k}=\left[\begin{array}{cc}0&A_2\\A_2&0\end{array}\right]$, & $_{\bf k}\gamma=\left[\begin{array}{cc}0&-B_2\\B_2&0\end{array}\right]$.
\end{tabular}$$

If we write
$$X_{\bf i}=\left[\begin{array}{cc}0&0\\0&\gamma_{\bf i}\end{array}\right],\hspace{5mm} _{\bf i}X=\left[\begin{array}{cc}0&0\\0& _{\bf i}\gamma\end{array}\right],$$
and so on, we get the following relations:
$$[X_{\bf i},S]= S{\bf i}\in\mathbb{H},\,\,\,[X_{\bf j},S]= S{\bf j}\in\mathbb{H},\,\,\,[X_{\bf k},S]=S{\bf k}\in\mathbb{H},$$
$$[ _{\bf i}X,S]= {\bf i}S\in\mathbb{H},\,\,\,[ _{\bf j}X,S]= {\bf j}S\in\mathbb{H},\,\,\,[ _{\bf k}X,S]= {\bf k}S\in\mathbb{H},$$
where $S{\bf i}$, $S{\bf j}$, $S{\bf k}$, ${\bf i}S$, ${\bf j}S$ and ${\bf k}S$ are well defined from $(\ref{ident})$.
Finally, we have that $\langle \gamma_{\bf i},\gamma_{\bf j},\gamma_{\bf k}\rangle$ and $\langle _{\bf i}\gamma, _{\bf j}\gamma, _{\bf k}\gamma\rangle$ are ideals of $\so(4)$, each one isomorphic to $\su(2)$. Thus, the adjoint representation of these components corresponds to right and left multiplication by immaginary quaternions, as claimed.
\end{profe}

Note that the construction in the above proof give us a direct way to verify the isomorphism $\so(4)=\su(2)\oplus\su(2)$.

In the remaining of this section we study the gradient vector fields $X_{(q,0,0)}$. To get a description about the bahavior of these vector fields, we start by studying the singularities of the vector fields $X_{(q,z,w)}$ on the sphere $S^3$. Remember that we can consider $S^3$ as the Lie group ${\rm SU}(2)$. So, given a vector field on ${\rm SU}(2)$ we define the function
$$F:{\rm SU}(2)\longrightarrow \su(2)$$
by setting
$$F(p)=d(D_{p^{-1}})_p(X(p))=X(p)\cdot p^{-1} . $$

Note that a point $p\in{\rm SU}(2)$ is a singular point for $X$ if, and only if, $F(p)=0$. Thus the set of singular points of the vector field $X$ is $F^{-1}\{0\}$ and $X$ has no singular points if, and only if, $F^{-1}\{0\}=\emptyset$, that is, $0$ does not belong to the image of the function $F$.

Now we prove the following lemma.
\begin{lema}

\item[(i)] The vector field $X$ is right-invariant if, and only if, its corresponding function $F$ is constant. More precisely, $F(p)=X(1)$ for all $p\in{\rm SU}(2)$.
\item[(ii)] The vector field $X$ is left-invariant if, and only if, $F(p)={\rm Ad}(p)X(1)$.

\end{lema}
\begin{profe}
In fact, for the item (i), suppose that $X$ is a right-invariant vector field, that is, 
$$d(D_g)_h(X(h))=X(hg),\hspace{5mm}\forall g,h\in{\rm SU}(2).$$
This implies that
$$F(p)=d(D_{p^{-1}})_p(X(p))=X(pp^{-1})=X(1)\in{\rm SU}(2).$$
By the other hand, if $F(p)=X(1)$ for all $p\in{\rm SU}(2)$, we get
$$F(p)=X(p)\cdot p^{-1}=X(1)\Longleftrightarrow X(p)=X(1)\cdot p,$$ 
that is, $X$ is right-invariant.

For item (ii), suppose first that $X$ is left-invariant, that is,
$$d(E_g)_h(X(h))=X(gh),\hspace{5mm}\forall g,h\in{\rm SU}(2).$$
Note that
$${\rm Ad}(p)X(1)=d(D_{p^{-1}})_p\circ d(E_p)_1(X(1))=d(D_{p^{-1}})_p(X(p))=F(p).$$
Reciprocally, if $F(p)={\rm Ad}(p)X(1)$, we have
$$F(p)=X(p)\cdot p^{-1}={\rm Ad}(p)X(1)=p\cdot X(1)\cdot p^{-1}\Longleftrightarrow X(p)=p\cdot X(1),$$
proving that $X$ is left-invariant and concluding the proof.    
\end{profe}

Now we take a look at the $F$-functions corresponding to the vector fields given in the previous theorem. A vector field 
$$X_{(q,z,w)}(p)=\frac{1}{2}(q-p\overline{q}p)+zp+pw,\hspace{5mm}p\in{\rm SU}(2),$$
can be written as 
$$X_{(q,z,w)}=X_{(q,0,0)}+X_{(0,z,0)}+X_{(0,0,w)}$$
and the corresponding function is given by $F_{(q,z,w)}=F_{(q,0,0)}+F_{(0,z,0)}+F_{(0,0,w)},$ since it is defined by the differential of the right translation by $p^{-1}$, which is linear. The functions $F_{(q,0,0)}$, $F_{(0,z,0)}$ and $F_{(0,0,w)}$ are given by
\begin{enumerate}
\item $X_{(0,z,0)}(p)=zp$ implies that $F_{(0,z,0)}(p)=X_{(0,z,0)}(p)\cdot p^{-1}=zpp^{-1}=z$, in other words, $F_{(0,z,0)}$ is constant. Note that this agrees with the item (i) above, since $X_{(0,z,0)}$ is a right-invariant vector field. 
\item Since $X_{(0,0,w)}=pw$ we have $F_{(0,0,w)}(p)=pwp^{-1}={\rm Ad}(p)X(1)$. And again it agrees with item (ii) above, because $X_{(0,0,w)}$ is left-invariant. Furthermore, as $p^{-1}=\overline{p}$, we get $F_{(0,0,w)}(p)=pw\overline{p}$.
\item Finally, for the vector field $X_{(q,0,0)}=\frac{1}{2}(q-p\overline{q}p)$ we have 
$$F_{(q,0,0)}(p)=X_{(q,0,0)}(p)\cdot p^{-1}=\frac{1}{2}(q-p\overline{q}p)\overline{p},$$
and so
$$F_{(q,0,0)}(p)=\frac{1}{2}(q\overline{p}-p\overline{q}p\overline{p})=\frac{1}{2}(q\overline{p}-p\overline{q})={\rm Im}q\overline{p}={\rm Im}p\overline{q}.$$
\end{enumerate}

It turns out that $F_{(q,0,0)}(p)={\rm Im}p\overline{q}=0$ if, and only if, $p\overline{q}=x\in\R$, that is, $p|q|^2=xq$, which means that the only singularities for $X_{(q,0,0)}$ are the antipodal points $p=\pm\frac{q}{|q|}$.

In the sequel we deal with the image of $F_{(1,0,0)}$. Since $F_{(1,0,0)}(p)={\rm Im}(p)$, we have that $F_{(1,0,0)}(S^3)$ is exactly the unit ball in ${\rm Im}\mathbb{H}$. To get a more detailed description of this image we look at the great circles passing by the elements $\pm{\bf 1}\in S^3$. The great circle $C_z$ is just the intersection of the plane generated by $\{{\bf 1},z\}$ with $S^3$ (here $z$ is an purely immaginary unit quaternion). 
The computations with $C_{{\bf i}}$ give us a good idea of how the general case works. Let $C_{{\bf i}}^+$ be the semicircle $C_{{\bf i}}^+=\{p\in C_{{\bf i}}\,|\,\langle p,i\rangle\geq0\}$. As $p$ runs through $C_{{\bf i}}^+$ from ${\bf 1}$ to $-{\bf 1}$ the values $F_{(1,0,0)}(p)$ cover twice the line segment $[0,{\bf i}]=\{t{\bf i}\,|\,t\in[0,1]\}$. In fact, for $p$ going from ${\bf 1}$ to ${\bf i}$ we have ${\rm Re}(p)\geq0$ and ${\rm Im}(p)$ goes from $0$ to ${\bf i}$. When $p$ goes from ${\bf i}$ to $-{\bf 1}$ we have ${\rm Re}(p)\leq0$ and ${\rm Im}(p)$ goes from ${\bf i}$ to $0$. Briefly, as $p$ goes from ${\bf 1}$ to $-{\bf 1}$ on $C_{{\bf i}}^+$ the image $F_{(1,0,0)}(p)$ goes from $0$ to ${\bf i}$ and then goes back from ${\bf i}$ to $0$, always lying on the line segment $[0,{\bf i}]$. On the other half $C^-_{{\bf i}}$ the situation is quite analogous, in this case however the image is the line segment $[-{\bf i},0]=\{-t{\bf i}\,|\,t\in[0,1]\}$, also covered twice as $p$ goes from ${\bf 1}$ to $-{\bf 1}$. By placing the things together we get
$$F_{(1,0,0)}(C_{{\bf i}})=[-{\bf i},{\bf i}]=\{t{\bf i}\,|\,t\in[-1,1]\}.$$   

An analogous reasoning works for a general great circle $C_z$, that is, $F_{(1,0,0)}(C_z)=[-z,z]$, the line segment joining $-z$ and $z$.

\section{Invariant control sets}

In this section we describe the control sets for the family of control systems on the sphere $S^{3}$ having 
$X_{(1,0,0)}$ as drift and control vector fields corresponding to pure quaternions. Others control systems, that is, considering others drifts and control vector fields the computations clearing become very intricate and the approach should be different to try determine the control sets.   The family of systems considered here is given by:
\begin{itemize}
\item[(i)] $X_{(1,0,0)}+uX_{(z,0,0)}$, where $z\in{\rm Im}\mathbb{H}$, $u\in[-1,1]$.
\item[(i')] $X_{(1,0,0)}+uX_{(z,0,0)}$, $z\in{\rm Im}\mathbb{H}$, $u=\pm1.$
\item[(ii)] $X_{(1,0,0)}+uX_{(z_1,0,0)}+vX_{(z_2,0,0)}$, $z_1,z_2\in{\rm Im}\mathbb{H}$ and $(u,v)\in[0,1]^2$.
\item[(ii')] $X_{(1,0,0)}+uX_{(z_1,0,0)}+vX_{(z_2,0,0)}$, $z_1,z_2\in{\rm Im}\mathbb{H}$, $(u,v)=(0,0),(1,0),(0,1)$.
\item[(iii)] $X_{(1,0,0)}+uX_{({\bf i},0,0)}+vX_{({\bf j},0,0)}+wX_{({\bf k},0,0)}$, $(u,v,w)\in B[0,\frac{1}{2}]\subset\R^3$.
\end{itemize}

Note that the Lie algebra $\so(1,4)$ cannot be generated by less than four symmetric elements, that is, elements in $\ms$. 

Hence we do not have controllability on ${\rm SO}(1,4)$ for these control systems whose vector fields correspond to elements in the component $\ms$ of the Cartan decomposition of $\so(1,4)$, that is, vector fields of the form $X_{(q,0,0)}$, $q\in\mathbb{H}$.

Now we introduce some necessary notations. A subset $C\subset S^3$ is said to be spherically convex when for any pair of points $p,q\in C$ every minimal geodesic segment joining them are contained in $C$. 
Given a subset $S\subset S^3$, denote by $K_S$ the cone spanned by $S$, that is, 
$$K_S=\{tp\,|\,p\in S,\,\,t\geq0\}\subset\ms,$$
and denote by ${\rm co}S$ the conic hull of $S$,
$${\rm co}S=\left\{\left.\sum_{i=1}^{n} a_ip_i\,\right|\,a_i\geq0,\,p_i\in S,\,n\in\N\right\}.$$ 
Cones are very useful in characterizing convex sets in the sphere, since it is well known that a non-empty subset $C\subset S^3$ is convex if, and only if, the cone $K_C$ is convex and pointed (that is, $K_C\cap(-K_C)\subseteq\{0\}$). This means that the proper convex sets on $S^3$ are the intersections of $S^3$ with pointed convex cones. In fact, if $K$ is a convex and pointed cone, then $K_C=K\cap S^3$ is convex since $K_C$ is exactly $K$. This means that the spherical convex hull of a subset $A\in S^3$ is the intersection of its conical hull ${\rm co}A$ with $S^3$.

With respect to the quaternions $\mathbb{H}$, the following notations are very useful. Given $p\in\mathbb{H}$, we write $p=p_0+p_1{\bf i}+p_2{\bf j}+p_3{\bf k}=(p_0,p_1,p_2,p_3)$ and we denote the imaginary part of $p$ by ${\vec p}=p_1{\bf i}+p_2{\bf j}+p_3{\bf k}$. Under these assumptions, $p\in\mathbb{H}$ can be seen as the sum of a scalar with a vector, that is, $p=p_0+{\vec p}$. Also, the product of two quaternions $p=p_0+p_1{\bf i}+p_2{\bf j}+p_3{\bf k}$ and $q=q_0+q_1{\bf i}+q_2{\bf j}+q_3{\bf k}$ can be written as 
\begin{equation*}
p\cdot q=p_0q_0-{\vec p}\cdot{\vec q}+p_0{\vec q}+q_0{\vec p}+{\vec p}\times {\vec q},\label{prod}
\end{equation*}
where $\vec{p}\cdot\vec{q}$ and $\vec{p}\times \vec{q}$ are identified with the inner and vector products of three-dimensional vectors.

Our first result is this section determines which spherical domes are invariant by $X_{(1\pm z,0,0)}$.

\begin{prop} Given $z\in{\rm Im}\mathbb{H}$, let $r_1=\sfrac{1}{\sqrt{1+|z|^2}}$. Then, for every $a\in[0,r_1)$ the spherical dome 
$$S_a=\{p\in S^3\,|\,{\rm Re}(p)\geq a\}$$
 is invariant under the family of vector fields $\{X_{(1+z,0,0)},X_{(1-z,0,0)}\}$.
\label{calotainv}
\end{prop}

\begin{profe}
For a pure quaternion $z$
$$X_{(1-z,0,0)}(p)=\frac{1}{2}\left(1-z-p^2-pzp\right).$$
If $p\in S^3\cap{\rm Im}\mathbb{H}$ the above expression implies
$$X_{(1-z,0,0)}(p)=\frac{1}{2}\left(2-z-pzp\right) \mbox{ and } X_{(1+z,0,0)}(p)=\frac{1}{2}\left(2+z+pzp\right).$$  
 
On the other hand, 
$pz=-{\vec p}\cdot{\vec z}+{\vec p}\times{\vec z}$ 
and hence $pzp=-({\vec p}\cdot{\vec z}){\vec p}+({\vec p}\times{\vec z})\times{\vec p}\in{\rm Im}\mathbb{H}$.


Since both $z$ and $pzp$ are pure quaternions, we can see that $X_{(1-z,0,0)}(p)$ and $X_{(1+z,0,0)}(p)$ have positive real parts (being $p$ a pure quaternion). This fact ensures that the upper half of the sphere $S^3$ is invariant under the trajectories of the family
$$\left\{X_{(1-z,0,0)},X_{(1+z,0,0)}\right\}.$$

Now, let $S^3\ni p=a+w$, with $0<a<1$ and $w\in{\rm Im}\mathbb{H}$. Note that $|w|=\sqrt{1-a^2}$. Then
\begin{eqnarray*}
X_{(1+z,0,0)}(p)&=&X_{(1+z,0,0)}(a+w)\\
&=&\frac{1}{2}\left(1+z-(a+w)^2+a^2z+azw+awz+wzw\right)\\
&=&{\rm Re}\left(X_{(1+z,0,0)}(a+w)\right)+{\rm Im}\left(X_{(1+z,0,0)}(a+w)\right),
\end{eqnarray*}
where the first and second part of the last sum are $\frac{1}{2}\left(1-a^2+{\vec w}\cdot{\vec w}-2a{\vec z}\cdot{\vec w}\right)$ and $\frac{1}{2}\left(z-2a{\vec w}-{\vec w}\times{\vec w}+a^2z+wzw\right)$ respectively. 

Similarly,
$$X_{(1-z,0,0)}(a+w)={\rm Re}\left(X_{(1-z,0,0)}(a+w)\right)+{\rm Im}\left(X_{(1-z,0,0)}(a+w)\right),$$
where the previous real and imaginary parts are respectively 
$$\frac{1}{2}\left(1-a^2+{\vec w}\cdot{\vec w}+2a{\vec z}\cdot{\vec w}\right) \mbox{ and } \frac{1}{2}\left(-z-2a{\vec w}-{\vec w}\times{\vec w}-a^2z-wzw\right).$$

The real parts of $X_{(1+z,0,0)}(a+w)$ and $X_{(1-z,0,0)}(a+w)$ are thus given by $1-a^2-a\sqrt{(1-a^2)}|z|\cos(t)$ and $1-a^2+a\sqrt{(1-a^2)}|z|\cos(t)$
where $t$ denotes the angle between $w$ and $z$.

Being $t$ the angle formed by $w$ and $z$, define the number 
$r_t=\frac{1}{\sqrt{1+|z|^2\cos^2(t)}}$.

Note that ${\rm Re}\left(X_{(1-z,0,0)}(a+w)\right)=0$ iff $a=r_t$ or $a=1$. For a fixed $t\in[0,2\pi]$ with $t\neq{\pi}/{2}$ and $t\neq{3\pi}/{2}$, we have 
\begin{eqnarray*}
{\rm Re}\left(X_{(1-z,0,0)}(a+w)\right)>0&\hbox{ if }&a\in[0,r_t)\\
{\rm Re}\left(X_{(1-z,0,0)}(a+w)\right)<0&\hbox{ if }&a\in\left(r_t,1\right) .
\end{eqnarray*} 
If $t={\pi}/{2}$ or $t={3\pi}/{2}$, then 
${\rm Re}\left(X_{(1-z,0,0)}(a+w)\right)>0$
for every $a\in[0,1)$, since in this case 
${\rm Re}\left(X_{(1-z,0,0)}(a+w)\right)=1-a^2.$

For a given $z\in{\rm Im}\mathbb{H}$ the number $r_t$ is minimum when $\cos(t)=1$. Observe that 
$r_1=\displaystyle\frac{1}{\sqrt{1+|z|^2}}$
 is precisely the real part of the singularities of $X_{(1\pm z,0,0)}$. This implies that if $a\in[0,r_1)$ then ${\rm Re}\left(X_{(1\pm z,0,0)}(a+w)\right)>0$ for every possible $w$. In other words, for every $a\in[0,r_1)$ the spherical dome $S_a=\{p\in S^3\,|\,{\rm Re}(p)\geq a\}$ is invariant under $\{X_{(1\pm z,0,0)}\}$.
\end{profe}

Now  we can describe the control sets for the families of control systems given in the beginning of this section. 

\textit{Case (i):} 

In this case note that
$$\Gamma=\left\{X_{(1,0,0)}+uX_{(z,0,0)}\,|\,z\in{\rm Im}\mathbb{H},\,\,u\in[-1,1]\right\},$$
since all the attractors for such a system lie on (and cover) the minimal geodesic segment $C$ in $S^3$ joining the attractors of $X_{(1+z,0,0)}$ and $X_{(1-z,0,0)}$, then  Theorem \ref{geral} implies that $C$ is the only invariant control set for $\Gamma$.

\textit{Case (i'):}

Using a geometrical approach we will show  that the invariant control set for the bang-bang control system 
$\{X_{(1,0,0)}\pm X_{(z,0,0)} , z\in{\rm Im}\mathbb{H} \}$
is  a geodesic segment joining $p_1$ and $p_2$, which are the singularities (attractors) of the vector fields $X_{(1+z,0,0)}$ and $X_{(1-z,0,0)}$, respectively. That is, the spherical closed convex hull $C$ of the set $\{p_1,p_2\}$, which is given by 
$$C={\rm cl}({\rm co}\{p_1,p_2\})\cap S^3.$$

\begin{teo}
On the sphere $S^3$ the invariant control set for the family $\{X_{(1\pm z,0,0)}\}$ is the minimal geodesic segment $C$ joining 
$$p_1=\frac{1+z}{|1+z|}\hbox{ and }p_2=\frac{1-z}{|1-z|}.$$
\end{teo}
\begin{profe}
Denote by $D$ the unique invariant control set for this system (it exists since $S^3$ is a compact Hausdorff space). Our aim is to show that $C=D$.

First note that for every pair of points $x,y\in C \setminus \{  p_1,p_2 \}$ we have a trajectory of the system joining $x$ and $y$. In fact, the solutions for this control system follow geodesic curves in $S^3$. If the solution with initial value $x$ that goes towards $p_1$ do not contains $y$, then the solution going towards $p_2$ does. Furthermore, $p_1\in{\rm cl}\mathcal{O}^+(p_2)$ and $p_2\in{\rm cl}\mathcal{O}^+(p_1)$. In other words, $C$ is contained in the closure of the positive orbit of all its elements. Also, the positive orbit of all its points is contained in $C$.

Now let $x\in D$. For any $y\in C$ we have $x$ in the closure of the positive orbit of $y$, that is, we can approximate $x$ by an trajectory of the control system starting at $y$. This trajectory is the concatenation of solutions, all of them lying on the great circle that contains ${\bf 1}$, $p_1$ and $p_2$, that is, all the trajectory is contained in this great circle. Since $C$ is spherically convex, every path of the trajectory is a minimal geodesic segment between its initial and final points. This implies that $x\in{\rm cl}(C)=C$, and the maximality of $D$ ensures that $C=D$. 
\end{profe}

\textit{Case (ii):}

Again, Theorem \ref{geral} implies that the invariant control set of the control system 
$\{X_{(1,0,0)}+uX_{(z_1,0,0)}+vX_{(z_2,0,0)} ;u,v\in[0,1] , z_1,z_2\in{\rm Im}\mathbb{H}\}$
is  the spherical closed convex hull of $A=\{{\bf 1},p_1,p_2,p_3\}$, where $p_1$ is the attractor corresponding to $u=1$ and $v=0$, $p_2$ corresponding to $u=0$ and $v=1$ and $p_3$ to $u=1$ and $v=1$. However, if the controls are in the interval $[-1,1]$ we need to add, in the set $A$, also the attractors corresponding to the control values $u,v=-1$. Any other attractors would belong to the spherical closed convex hull too, and a similar reasoning works.

\textit{Case (ii')}:

In this case with only two vector fields for the control, the situation is quite similar to that one in case (i') and a geometrical reasoning still works. 

\begin{teo}
If $\Gamma=\{X_{(1,0,0)},X_{(1+z_1,0,0)},X_{(1+z_2,0,0)}\}$, $z_1,z_2\in{\rm Im}\mathbb{H}$, is a family of vector fields on the sphere $S^3$ such that the attractors $p_1$ and $p_2$ for $X_{(1+z_1,0,0)}$ and $X_{(1+z_2,0,0)}$ do not belong to the same great circle containing ${\bf 1}$, then the invariant control set for $\Gamma$ in $S^3$ is the spherical closed convex hull of $A=\{{\bf 1},p_1,p_2\}$.
\end{teo}

\begin{profe}
Let ${\rm co}A$ be the conic hull of $A$. The closed convex hull of $A$ is 
$$C={S^3\cap{\rm cl}({\rm co}A)}.$$
Note that $C$ has empty interior in $S^3$ since it is a two-dimensional submanifold (${\rm co}A$ is a cone generated by three elements). To avoid confusion, write ${\rm algint}C$ for the algebraic interior of $C$. The relative boundary of ${\rm co}A$ is the union of the conic hulls of $\{{\bf 1},p_1\}$, $\{{\bf 1},p_2\}$ and $\{p_1,p_2\}$, that is, the relative boundary of $C$ is union of great circle arcs joining ${\bf 1}$ to $p_1$, ${\bf 1}$ to $p_2$ and $p_1$ to $p_2$.

To see that $C$ is the desired invariant control set, let $x\in{\rm algint}C$. Denote by $X_1(x)$, $X_{p_1}(x)$ and $X_{p_2}(x)$ the solutions for the corresponding vector fields with initial value $x$. Remember that all these solutions follows geodesics on $S^3$. Let $C_1$ be the closed convex hull of $\{x,{\bf 1},p_1\}$, $C_2$ the closed convex hull of $\{x,{\bf 1},p_2\}$
and $C_3$ that one of $\{x,p_1,p_2\}$. It is clear that $C=C_1\cup C_2\cup C_3$. If $y\in{\rm algint}C_1$, then we obtain a trajectory from $x$ to $y$ in the following way. The convexity of $C$ tells us that the solution $X_1(y)$ meets $X_{p_1}(x)$ in $y_1\in{\rm algint}C$ (the intersection $X_1(y)\cap X_{p_1}(x)$ is nonvoid since both solutions are contained in the same three dimensional subspace of $\ms$, and the intersection belongs to ${\rm algint}C$ by convexity). The solution $X_1(y_1)$ contains $y$ since $X_1(y_1)$ and $X_1(y)$ are contained in the same great circle of $S^3$. Thus, joining $X_{p_1}(x)$ and $X_1(y_1)$ we get a trajectory from $x$ to $y$. Analogous approach work for ${\rm algint}C_2$ and ${\rm algint}C_3$, just interchanging the considered solutions. If $y\in{\rm algint}C$ does not belong to ${\rm algint}C_i$, $i=1,2,3$, then $y$ belongs to $X_1(x)$, $X_{p_1}(x)$ or $X_{p_2}(x)$, and can be trivially reached from $x$. This proves that $C$ is contained in the closure of the positive orbit of all its elements. 
 
Now, let $D$ be the unique invariant control set for this control system and take $x\in D$. We just proved that $C\subset D$, by the maximality of $D$. That is, for any $y\in C$ we have $x\in{\rm cl}\mathcal{O}^+(y)$. In other words, we can get close of $x$ by trajectories starting at $y$. Every path of such an trajectory is a minimal geodesic segment, and this implies that $x\in{\rm cl}C$, by the very definition of spherical closed convex hull. Hence $x\in C$, and so $C=D$.
\end{profe}

\textit{Case (iii):}

Note that the control systems considered above does not satisfy the LARC, as said in the beginning of this section. One can note that all the invariant control sets obtained have empty interior in $S^3$. The next theorem present a study of a control system satisfying the LARC. Note that checking LARC is easier if it is done with matrices instead of vector fields (see Proposition \ref{infiniact}).
\begin{teo}
Consider the control system 
$$X_{(1,0,0)}+uX_{({\bf i},0,0)}+vX_{({\bf j},0,0)}+wX_{({\bf k},0,0)},\hspace{5mm}(u,v,w)\in B[0,\sfrac{1}{2}]\subset\R^3,$$
where $B[0,\sfrac{1}{2}]$ stands for the closed ball of radius $\frac{1}{2}$ in $\R^3$ centered at the origin. Then, the invariant control set for this control system is exactly the spherical dome 
$$C=\{p\in S^3\,|\,{\rm Re}(p)\geq\sfrac{1}{\sqrt{2}}\}.$$
\end{teo}
\begin{profe} 
To prove this theorem, note first that the invariance follows from computations similar to those we have done for $X_{(1\pm z,0,0)}$ in Proposition \ref{calotainv}. The singularities corresponding to the vector fields $X_{(1\pm{\bf i},0,0)}$, $\pm X_{(1\pm{\bf j},0,0)}$ and $X_{(1\pm{\bf k},0,0)}$ on the invariant hemisphere are
$$\frac{1\pm{\bf i}}{\sqrt{2}},\hspace{2mm}\frac{1\pm{\bf j}}{\sqrt{2}},\hspace{2mm}\frac{1\pm{\bf k}}{\sqrt{2}}$$
and such singularities lie on the boundary of the spherical dome $C$, since their real parts are all $\sfrac{1}{\sqrt{2}}$, and the same statement is true for the singularities corresponding to vector fields
$$X_{(1,0,0)}+uX_{({\bf i},0,0)}+vX_{({\bf j},0,0)}+wX_{({\bf k},0,0)}$$
such that $u^2+v^2+w^2=1$. 

By the other hand, let $p$ in the boundary of $C$. It must have the form $p=\sfrac{1}{\sqrt{2}}+z$, with $z\in{\rm Im}\mathbb{H}$. Observe that
$$|p|^2=\frac{1}{2}+|z|^2\Longrightarrow|z|=\frac{1}{\sqrt{2}}.$$
Write $z=u_0{\bf i}+v_0{\bf j}+w_0{\bf k}$. Then $u_0^2+v_0^2+w_0^2=\sfrac{1}{2}$ and the attractor corresponding to the vector field
$$X_{(1,0,0)}+\sqrt{2}u_0X_{({\bf i},0,0)}+\sqrt{2}v_0X_{({\bf j},0,0)}+\sqrt{2}w_0X_{({\bf k},0,0)}$$
is the point $p$, as can be easily checked. This proves that every point belonging to the boundary of $C$ is a singularity of a vector field of the control system corresponding to a control $(u,v,w)$ such that $u^2+v^2+w^2=1$.

The previous reasoning still holds for controls $(u,v,w)\in B[0,a]$, with $0<a<1$. In this case, the singularities corresponding to controls $(u,v,w)$ such that $u^2+v^2+w^2=a^2$ are in the boundary of the spherical dome 
$$C_a=\left\{p\in S^3\,\left|\,{\rm Re}(p)\geq\frac{1}{\sqrt{1+a^2}}\right.\right\},$$
and so we conclude that all of the points belonging to $C$ are singularities corresponding to some vector field of the control system. The result is thus a consequence of Theorem \ref{geral}.
\end{profe}


\end{document}